\documentclass[a4paper,10pt]{amsart}
\usepackage[english]{certus}
\usepackage{a4wide}
\usepackage{amsmath,amsfonts,amssymb,amsthm}
\usepackage{comment}
\usepackage{tikz-cd}
\usepackage{bbold}

\theoremstyle{definition}

\newtheorem*{As*}{Assumption}
\theoremstyle{plain}
\newtheorem*{Claim*}{Claim}
\theoremstyle{remark}
\DeclareMathOperator{\spa}{span}

\DeclareMathOperator{\SL}{SL}
\DeclareMathOperator{\PSL}{PSL}

\newcommand{\NN}{\ensuremath{\mathbb{N}}}
\newcommand{\ZZ}{\ensuremath{\mathbb{Z}}}
\newcommand{\RR}{\ensuremath{\mathbb{R}}}
\newcommand{\CC}{\ensuremath{\mathbb{C}}}

\title{A rigidity result for normalized  subfactors}
\author{Vadim Alekseev}
\address{Vadim Alekseev, Technische Universit\"{a}t Dresden, Fachrichtung Mathematik, Institut f\"{u}r Geometrie, 01062 Dresden, Deutschland}
\email{vadim.alekseev@tu-dresden.de}

\author{Rahel Brugger}
\address{Rahel Brugger, Technische Universit\"{a}t Dresden, Fachrichtung Mathematik, Institut f\"{u}r Geometrie, 01062 Dresden, Deutschland}
\email{rahel.brugger@tu-dresden.de}


\begin{document}
\onehalfspace

\begin{abstract}
		We  show a rigidity result for subfactors that are normalized by a representation of a lattice $\Gamma$ in a higher rank simple Lie group with trivial center into a finite factor. This implies that every subfactor of $L\Gamma$ which is normalized by the natural copy of $\Gamma$ is trivial or of finite index.
		This article is based on some results from the second named author's Ph.D. thesis.
\end{abstract}

\maketitle
\section{Introduction}
It seems like a natural generalization of Margulis' Normal Subgroup Theorem to ask
wheter every regular subfactor of the group von Neumann algebra of a lattice in a higher-rank simple Lie group with trivial center is trivial or of finite index.

In this article we make a small step into the direction of answering this question by looking at the special case where a subfactor $N\subset L\Gamma$ is actually normalized by a unitary representation $\pi\colon\Gamma\to U(L\Gamma)$ such that $N$ and $\pi(\Gamma)$ generate $M$.
We use methods developed by Jesse Peterson for the proof of his character rigidity theorem to prove the following theorem.
\begin{Thm}[Theorem \ref{Thm from Pet machine}]
	Let $\Gamma$ be a lattice in a simple real Lie group $G$ which has trivial center and real rank at least $2$. Let $M$ be a finite factor, $N\subset M$ a subfactor and $\pi\colon \Gamma \to \mathcal N_M(N)$ a unitary representation of $\Gamma$ into the normalizer of $N$ such that the action $\Gamma\curvearrowright M$ given by $\alpha_{\gamma}(x)=\pi(\gamma)x\pi(\gamma^{-1})$ is ergodic and $M=(N\cup\pi(\Gamma))''$.
	Then $M\ne N\rtimes \Gamma$ or  $[M:N]<\infty$.
\end{Thm}

Peterson's proof is inspired by Margulis' proof in the sense that the proof of the normal subgroup theorem is based on the fact that an amenable discrete group with property (T) is finite, whereas the proof of character rigidity is based on the fact that an amenable factor with property (T) is finite-dimensional.

We adjust Peterson's proof to the situation of subfactors described above by putting coefficients in $N$ into it.
Then we use that if an inclusion $N\subset M$ is both coamenable and corigid and the relative commutant is finite-dimensional, the inclusion is of finite index.

Another approach to the result above would be to study the character
$\gamma\mapsto \|E_N(\pi(\gamma))\|_2^2$. If this happens to be extremal, the above theorem follows directly from character rigidity \cite[Theorem C]{PetersonRigidity2016}: if the GNS construction of this character generates a finite dimensional von Neumann algebra, then the index is finite; if it is the regular character, then $M$ is the crossed product $N\rtimes\Gamma$. However, as we observe in the last section, no direct reduction to the extremal case seems possible: the above character happens to be extremal if and only if it is the regular character.

\subsection*{Acknlowledgements}
We would like to thank Sayan Das, Jesse Peterson and Andreas Thom for helpful comments on an earlier version of this text.
\section{Preliminaries}
Let $M$ for the rest of this text be a finite factor with trace $\tau$, $N\subset M$ a subfactor  and let $\Gamma$ be a discrete group. With $\mathcal{N}_M(N):=\{u\in U(M)| \,uNu^*=N\}$ we denote the normalizer of $N$ inside $M$. $[M:N]:=\dim_N(L^2M)$ is the index of $N\subset M$. We denote by $J$ the antilinear, bounded operator on $L^2(M)$ that extends $x\mapsto x^*$ for $x\in M$.
\subsection{Coamenability and Corigidity of inclusions of von Neumann algebras}
Amenability and property (T) of $M$ in Peterson's proof will be replaced by coamenability and corigidity of the inclusion $N\subset M$.
\begin{Def}[{\cite[3.2.3 (ii)]{Popa}}]
		Let $M$ be a factor and $N$ a von Neumann subalgebra. Then $N\subset M$  is \emph{coamenable} if there exists a conditional expectation $E\colon B(L^2 M) \cap (JNJ)' \to M$.
\end{Def}
\begin{Def}[\cite{Popa06}]
	Let $M$ be a factor and $N$ a von Neumann subalgebra. Then $N\subset M$ is  called \emph{corigid} if every $M$-bimodule $H$ with $N$-central norm one vectors $\xi_n\in H$  such that $\| x\xi_n - \xi_n x\|\rightarrow 0$ for all $x\in M$ contains a non-zero $M$-central vector.
\end{Def}
Note that  in \cite{Popa} coamenability is called amenability and corigidity is called rigidity.
\begin{Thm}[{\cite[4.1.8 (iv)]{Popa}}]\label{Popa}
	If an inclusion $N\subset M$ is coamenable and corigid and $N'\cap M$ is finite dimensional, then the inclusion is of finite index.
\end{Thm}
\subsection{Actions on von Neumann algebras}
\begin{Def}
	An action $\sigma\colon\Gamma\to\Aut(M)$ of a group on a von Neumann algebra is \textit{ergodic} if the fixed point algebra is $\CC$.
\end{Def}
Let us recall the definition und some properties of induced actions on von Neumann algebras.
\begin{Def}\label{induced representation}	
	Let $\Gamma\subset G$ be a closed subgroup of a locally compact group and $\theta\colon\Gamma\to\Aut(M)$ a continuous action. 
	Pick a Borel section $s\colon G/\Gamma \to G$ and let $\chi\colon G\times G/\Gamma \to \Gamma$ be the cocycle given by $\chi(g,x)=s(gx)^{-1}gs(x)$.
	
	Then the \emph{induced action}	$\tilde{\theta}$  of $G$ on $L^\infty(G / \Gamma) \overline{\otimes} M$, which we view as bounded functions from $G/\Gamma$ to $M$, is given by 	\begin{align*}
	\tilde{\theta}_g(f)(x) := \theta_{\chi(g,g^{-1}x)}f(g^{-1} x),	
	\end{align*}
	for $g\in G$, $f\in L^\infty(G / \Gamma) \overline{\otimes} M$ and $x\in G/\Gamma$.
\end{Def} 
\begin{Rem}\label{induced action other picture}
	Let $R$ be the $G$-action on $L^\infty(G)$ given by right multiplication. Then
	\[\Psi\colon L^\infty(G/\Gamma)\vntens M\to (L^\infty(G)\vntens M)^{(R\otimes\theta)(\Gamma)},\quad \Psi(f)(g)=\theta_{s(\Gamma)\chi(g,g\Gamma)}(f(g\Gamma))\]
	is an isomorphism and 
	\[\Psi(\tilde{\theta}_g(f))=L\otimes\id(g)\Psi(f),\]
	where $L$ is the $G$-action on $L^\infty(G)$ given by left multiplication.
\end{Rem}
\begin{Lemma}\label{induced representation ergodic}
$(L^\infty(G/\Gamma)\vntens M)^{\tilde{\theta}(G)}\cong1\otimes M^{\theta(\Gamma)}$. In particular, $\tilde{\theta}$ is ergodic iff $\theta$ is.
\end{Lemma}
\begin{proof}
	By Remark \ref{induced action other picture},
	\begin{align*}
	(L^\infty(G/\Gamma)\vntens M)^{\tilde{\theta}(G)}
	\cong (L^\infty(G)\vntens M)^{(R\otimes\theta)(\Gamma)\cup(L\otimes\id)(G)}=1\otimes M^{\theta(\Gamma)},
	\end{align*}
	hence $(L^\infty(G/\Gamma)\vntens M)^{\tilde{\theta}(G)}=\CC$	if and only if $M^{\theta(\Gamma)}=\CC$.
\end{proof}

\section{A question about regular subfactors of the von Neumann algebra of lattices in higher-rank groups}\sectionmark{A question}

We want to study possible analogues of Margulis' Normal Subgroup Theorem \cite[Theorem IX.5.3]{MargulisDiscrete1991} in the setting of subfactors. 
\begin{Thm}[Margulis' Normal Subgroup Theorem]
	Let $\Gamma$ be an irreducible lattice in a higher-rank simple Lie group $G$ with trivial center. Then every normal subgroup of $\Gamma$ is trivial or of finite index.
\end{Thm}
A typical example of such a group is $\PSL(n,\ZZ)\subset\PSL(n,\RR)$ for $n\geq 3$.
Margulis' Theorem was generalized by J. Peterson in \cite{PetersonRigidity2014} as follows.
\begin{Thm}[Peterson]
	Let $G$ be a property (T) semi-simple Lie group with trivial center, no compact factors, and real rank at leasdt 2, and let $\Gamma < G$ be an irreducible lattice in $G$.  Then for every unitary representation $\pi$ of $\Gamma$ such that $\pi(\Gamma)''$ is a finite factor $\pi$ extends to an isomorphism $L\Gamma\to \pi(\Gamma)''$ or  $\pi(\Gamma)''$ is finite-dimensional.
\end{Thm}
It should be possible to do everything in this article with the same assumptions on $\Gamma$ as in the above theorem.  We restrict ourselves to the simple real case to avoid some technicalities.

When replacing groups by factors the analogue of a normal subgroup is a regular subfactor.
\begin{Def}
	An inclusion of von Neumann algebras $N\subset M$ is \emph{regular} if the normalizer of $N$ generates $M$, i.e., $\mathcal{N}_M(N)''=M$.
\end{Def}
\begin{Qu}\label{subfactor rigidity}
	Is it true that if $\Gamma$ is as above and $N\subset L\Gamma$ a regular subfactor then $N=\CC$ or $[L\Gamma:N]<\infty$?
\end{Qu}
This question has probably been asked before, but we couldn't find it in the literature.

In the following we restrict our attention to the situation where the image of $\Gamma$ is not only in $\mathcal{N}_M(N)''$, but even in $\mathcal{N}_M(N)$ in order to make the question accessible to Peterson's methods from the proof of character rigidity.
We allow subfactors in von Neumann algebras a bit more general than $L\Gamma$.
\begin{As*}
	For the rest of this article let $\Gamma$ be a lattice in a simple real Lie group $G$ which has trivial center and real rank at least 2. Let $M$ be a finite factor with trace $\tau$, $N\subset M$ a subfactor and $\pi\colon \Gamma \to \mathcal N_M(N)$ a representation of $\Gamma$ into the normalizer of $N$ such that the action $\Gamma\curvearrowright M$ given by $\alpha_{\gamma}(x)=\pi(\gamma)x\pi(\gamma^{-1})$ is ergodic and $M=(N\cup\pi(\Gamma))''$.
\end{As*}

\begin{Ex}
	$M=L\Gamma$ with $\lambda:\Gamma\to L\Gamma$ the left regular representation and $N\subset L\Gamma$ a subfactor which is normalized by $\lambda(\Gamma)$ is as in the assumption. $M$ is a factor and the conjugation action is ergodic because $\Gamma$ is i.c.c.\@.
\end{Ex}

\section{Peterson machine with coefficients in $N$}
In this section we adjust the proof of Peterson's character rigidity theorem in \cite{PetersonRigidity2014} and \cite{PetersonRigidity2016} to the situation described above by putting coefficients in $N$ into it. Setting $N=\CC$ gives back the proof of character rigidity.

We will need a bunch of subgroups, which we first define in the case of $G=\SL (n,\RR)$.
\begin{Ex}
	For $G=\SL (n,\RR)$ let $P$ be the subgroup of upper triangular matrices and $V$ the subgroup of upper triangular matrices with 1 on the diagonal. Fix numbers $0=j_0<j_1<j_2<\dots<j_k=n$. We define  now subgroups consisting each of all block matrices in $\SL (n,\RR)$ of a certain structure:
	{\allowdisplaybreaks
	\begin{align*}
	P_0:=&\left\lbrace\begin{pmatrix}
	A_{11}& A_{12} & \dots & A_{1k}\\
	0& A_{22}&\dots & A_{2k}\\
	\vdots & \vdots & \ddots & \vdots\\
	0& 0& \dots & A_{kk}
	\end{pmatrix} \right\rbrace, \quad
	V_0:=\left\lbrace\begin{pmatrix}
	\mathbb{1} & A_{12} & \dots & A_{1k}\\
	0& \mathbb{1} &\dots & A_{2k}\\
	\vdots & \vdots & \ddots & \vdots\\
	0& 0& \dots & \mathbb{1}
	\end{pmatrix} \right\rbrace,  \\
	R_0:=&\left\lbrace\begin{pmatrix}
	A_{11}& 0 & \dots & 0\\
	0& A_{22}&\dots & 0\\
	\vdots & \vdots & \ddots & \vdots\\
	0& 0& \dots & A_{kk}
	\end{pmatrix} \right\rbrace,
	\quad
	L_0:=\left\lbrace\begin{pmatrix}
	V_{11}& 0 & \dots & 0\\
	0& V_{22}&\dots & 0\\
	\vdots & \vdots & \ddots & \vdots\\
	0& 0& \dots & V_{kk}
	\end{pmatrix} \right\rbrace,
	\end{align*}}
	where $A_{il}$ are arbitrary matrices of size $(j_i-j_{i-1})\times (j_l-j_{l-1})$, $V_{ii}$ are upper triangular matrices with $1$ on the diagonal and $\mathbb{1}$ is an identity matrix of fitting size.
	
	For each of these subgroups we denote by $\overline{P}, \overline{V}$, ect. the corresponding transposed subgroup.
\end{Ex}
\begin{Def}
	For general $G$, let $S$ be an $\RR$-split maximal torus, $P$ a minimal parabolic subgroup containing $S$ and $V<P$ its unipotent radical. Let $\overline{P}$ be the opposite parabolic and $\overline{V}$ its unipotent radical. Let $P_0$ be another parabolic subgroup s.t. $P<P_0\lneq G$. Let $V_0$ be the unipotent radical of $P_0$ and $\overline{P_0}, \overline{V_0}$ the corresponding opposite subgroups.
	Let $R_0$ be the reductive component of $P_0$ containing $S$ so that $P_0=R_0\rtimes V_0$ and $\overline{L_0}=R_0\cap\overline{V}$. Then $\overline{V}=\overline{V_0}\rtimes\overline{L_0}$. See \cite[I.1.2]{MargulisDiscrete1991} for the definitions.
\end{Def}

We have the following commuting diagram:
\[
\begin{tikzcd}
\overline{V}= \overline{V_0}\rtimes\overline{L_0}\arrow{r}{(v,l)\mapsto v} \arrow[swap]{d}{v\mapsto vP} & \overline{V_0} \arrow{d}{v\mapsto vP_0} \\
G/P  \arrow{r}{gP\mapsto gP_0} & G/P_0.
\end{tikzcd}
\]
\begin{Lemma}
	The vertical maps map measures in the class of the Haar measure to $G$-quasiinvariant measures. They are measure isomorphisms when equipping the quotient spaces with these measures.
\end{Lemma}
\begin{proof}
	Let $\mu_{\overline V}$ be a left Haar measure on $\overline V$ and let $\lambda\in\mathcal{M}(G/P)$ be the image of $\mu_{\overline{V}}$. This  defines by \cite[Proposition VII.2.1.4]{BourbakiIntegration2004} a measure $\lambda^\#$ on $G$ given by
	\[\int_G f d\lambda^\#=\int_{G/P}\int_P f(gp)d\mu_P(p)d\lambda(gP),
	\]
	where $\mu_P$ is a left Haar measure on $P$. By \cite[Lemma IV.2.2]{MargulisDiscrete1991}, the map $(\overline{V}\times P, \mu_{\overline{V}}\otimes\mu_P)\to (G,\mu_G),\, (v,p)\mapsto vp^{-1}$, is a homeomorphism onto the image and a measure isomorphism, $\mu_G$ being a suitably normalized left Haar measure on $G$. This implies that $\lambda ^\#=( 1_{\overline{V}}\otimes \Delta _P) \cdot \mu_G$, where $\Delta_P$ is the modular function on $P$. It follows by \cite[Lemma VII.2.5.4]{BourbakiIntegration2004} that $\lambda$ is $G$-quasiinvariant.
	Now the same follows for the images of measures that are strongly equivalent to Haar measure and analogously for such measures on $\overline{V_0}$.
	
	The vertical maps are then measure isomorphisms because they map the measures to each other and are injective.
\end{proof}
Let $\nu$ and $\rho$ be probability measures on $\overline{V_0}$ resp. $\overline {L_0}$ in the class of the Haar measure; the image of $\nu$ on $G/P_0$ is still denoted by $\nu$. We equip $\overline V = \overline{V_0}\rtimes \overline{L_0}$ with the product measure $\nu\times \rho$.

Let $G$ act on $\overline{V}$ and $\overline{V_0}$ in the way that makes the above diagram $G$-equivariant. This transforms the action of $\overline{V}$ on $G/P$ to left multiplication on $\overline{V}$ and the action of $R_0$ on $G/P$ to the action induced by conjugation on $\overline{V}$.

Let $\sigma$ be the corresponding action of $\Gamma$ on $L^\infty(G/P)$ and $\sigma^0$ the corresponding Koopman representation on $L^2(G/P)$. Let
\[P_1:=1\otimes P_{\hat{1}}\in L^\infty(G/P)\overline{\otimes} B(L^2M),\]
where $P_{\hat{1}}$ is the orthogonal projection on $\CC\hat{1}\subset L^2M$ with $M$ as in the assumption. 
Let
\[\mathcal{B}:=(L^\infty(G/P)\overline{\otimes} B(L^2M))\cap\{\sigma_\gamma\otimes(J\pi(\gamma)J)|\gamma\in\Gamma\}'\cap(1\otimes JNJ)'.\]
\begin{Lemma}\label{cond exp on B}
	There exists a conditional expectation 
	\[E\colon (L^\infty(G/P)\overline{\otimes} B(L^2M))\cap(1\otimes JNJ)'\to \mathcal{B}.\]
\end{Lemma}
\begin{proof}
	Let $H=L^2(M)$ and let $\theta\colon \Gamma \to \Aut(B(H))$ be conjugation with $J\pi(\cdot)J$.
	Let $\tilde{\theta}$ be the induced action of $G$ on $L^\infty(G / \Gamma) \overline{\otimes} B(H)$ as in Definition \ref{induced representation} with a section $s\colon G/\Gamma \to G$ and $\chi\colon G\times G/\Gamma \to \Gamma$ given by $\chi(g,x)=s(gx)^{-1}gs(x)$.
	$\tilde{\theta}$ is also well-defined on 
	$L^\infty(G / \Gamma)\overline{\otimes} B(H) \cap (1\otimes JNJ)'$,	
	which we view as bounded functions from $G/\Gamma$ to $B(H) \cap (JNJ)'$. To see this let $f\in L^\infty(G / \Gamma)\overline{\otimes}(B(H) \cap (JNJ)')$, $x\in G/\Gamma$, $n\in N$, $g\in G$, $\gamma= \chi(g,g^{-1}x)$ and calculate
	\begin{align*}
	\tilde{\theta}_g(f)(x) JnJ 
	&= J\pi(\gamma)J f(g^{-1}x) J\pi(\gamma^{-1})J JnJ \\
	&= J\pi(\gamma)J f(g^{-1}x) J\alpha_{\gamma^-1}(n)J J\pi(\gamma^{-1})J \\
	&= J\pi(\gamma)J	J\alpha_{\gamma^{-1}}(n)J f(g^{-1}x)  J\pi(\gamma^{-1})J \\
	&= JnJ J\pi(\gamma)J f(g^{-1}x)J\pi(\gamma^{-1})J \\
	&= JnJ \tilde{\theta}_g(f)(x).
	\end{align*}	
	
	$P$ is amenable \cite[IV.4.4]{MargulisDiscrete1991}, hence \cite[Theorem 7.4]{PetersonRigidity2016} gives that there is a conditional expectation
	\begin{align*}
	E\colon L^\infty(G /\Gamma)\otimes( B(H) \cap (JNJ)')
	\to\, &(L^\infty(G /\Gamma)\otimes( B(H) \cap (JNJ)'))^{\tilde{\theta}(P)} \\
	&= (L^\infty(G / \Gamma)\otimes B(H))^{\tilde{\theta}(P)} \cap (1\otimes JNJ)'.
	\end{align*}
	But 
	\begin{align*}
	(L^\infty(G / \Gamma)\overline{\otimes} B(H))^{\tilde{\theta}(P)}
	&\cong (L^\infty(G)\overline{\otimes} B(H))^{(L\otimes\id(P))\times (R\otimes\theta(\Gamma))}\\
	&\cong (L^\infty(G)\overline{\otimes} B(H))^{(R\otimes\id(P))\times (L\otimes\theta(\Gamma))}\\
	&\cong (L^\infty(G/P)\overline{\otimes} B(H))^{\sigma\otimes\theta(\Gamma)}.
	\end{align*}
	Here the first isomorphism is the map $\Psi$ given in Remark \ref{induced action other picture} and
	the second isomorphism is $f\mapsto (g\mapsto f(g^{-1}))$.	
	Thus \[(L^\infty(G / \Gamma)\overline{\otimes}	B(H))^{\tilde{\theta}(P)} \cap (JNJ)'
	=(L^\infty(G/P)\overline{\otimes} B(H))^{\sigma\otimes\theta(\Gamma)}\cap (1\otimes JNJ)'
	=\mathcal{B}.\]
\end{proof}

The following lemma is \cite[Lemma 4.4]{PetersonRigidity2014}, only with different $\mathcal{B}$, which does not change the proof. We give it anyway in order to provide more details.
\begin{Lemma}\label{lem 11.1}
	Let \[x=x^*\in \mathcal{B}\subset L^\infty(G/P)\overline{\otimes} B(L^2M)
	=L^\infty(\overline{V_0})\overline{\otimes}L^\infty(\overline{L_0})\overline{\otimes} B(L^2M)\] and view it as a function from $\overline{V_0}$ to $L^\infty(\overline{L_0})\overline{\otimes} B(L^2M)$.
	Let \[x_0\in L^\infty(\overline{L_0})\overline{\otimes} B(L^2M)\] be in the SOT-essential range of $x$. 
	Then there exists a $y=y^*\in\mathcal{B}$ such that $yP_1\in L^\infty(\overline{L_0})\overline{\otimes} B(L^2M)$ and 
	$P_1 yP_1=P_1x_0P_1$.
\end{Lemma}
\begin{proof}
	That $x_0$ is in the SOT-essential range of $x$ means that there are subsets $E_j\subset\overline{V_0}$ of positive measure such that for all $\eta\in L^2(M),\,\xi_L\in L^2(\overline{L_0})$ and $\epsilon>0$ there exists an $N$ with
	\begin{align}\label{essential range}
	\int_{\overline{L_0}}\|(x(v,l)-x_0(l))\eta\|^2|\xi_L(l)|^2\,d\rho(l) <\epsilon
	\end{align} 
	for all $j>N$ and all $v\in E_j$. By the proof of \cite[Lemma 4.3]{PetersonRigidity2014} there are $\gamma_j\in\Gamma$ and $h_j\in \overline{V_0}\rtimes\mathcal{Z}(R_0)$ such that $\gamma_jh_j^{-1}\to e$ and  $\nu(h_j E_j)\to 1$.
	We first show that these can be chosen in a way that $\sigma_{\gamma_j}(x)\to x_0$ in SOT. Take a countable SOT-basis of neighborhoods of zero in the unit ball of $L^\infty(\overline{V_0})\overline{\otimes} L^\infty(\overline{L_0})\overline{\otimes} B(L^2M)$, denoted by $\{U_j\}_{j\in\NN}$, such that $U_j\searrow \{0\}$. As the action $\sigma$ is strongly continuous, there are numbers $k(j)\in \mathbb N$ and neighborhoods $e\in O_j\subset G$ such that $O_j\subset O_i$ if $j>i$ and 
	\begin{align}\label{O_j}
	\sigma_g(x_0+U_{k(j)})\subset x_0+U_j\quad\forall g\in O_j.
	\end{align}
	We can choose the $\gamma_j$ in \cite[Lemma 4.3]{PetersonRigidity2014} in a way that $\gamma_jh_j^{-1}\in O_j$ for all $j$. We will show now first that $\sigma_{\gamma_j}(x)\to x_0$ if $\sigma_{h_j}(x)\to x_0$ and then that $\sigma_{h_j}(x)\to x_0$, all in SOT.
	
	So assume that $\sigma_{h_j}(x)\to x_0$, hence $\forall j\; \exists N:\sigma_{h_i}(x)\in x_0+ U_{k(j)}$ for all $ i>N $. Then by (\ref{O_j})
	$\sigma_{\gamma_i}(x_0)=\sigma_{\gamma_i h_i^{-1}}\sigma_{h_i}(x_0)\in x_0+ U_j$ for all $i>\max\{N,j\}$ because $\gamma_i h_i^{-1}\in O_j$. So then $\sigma_{\gamma_j}(x)\to x_0$.
	
	To show that $\sigma_{h_j}(x)\to x_0$ let $\eta\in L^2M,\xi_L\in L^2(\overline{L_0}),\xi_V\in L^2(\overline{V_0})$. Then, as the $h_j\in \overline{V_0}\rtimes\mathcal{Z}(R_0)$ act trivially on $\overline{L_0}$ and using (\ref{essential range}),
	\begin{align*}
	&\|1_{h_j E_j}(\sigma_{h_j}(x)-x_0)(\xi_V\otimes\xi_L\otimes\eta)\|^2\\
	=&\int_{h_jE_j\times\overline{L_0}}\|(\sigma_{h_j}(x)(v,l)-x_0(l))\eta\|^2|\xi_V(v)\xi_L(l)|^2\, d\nu(v)d\rho(l)\\
	=&\int_{h_jE_j\times\overline{L_0}}\|(x(h_j^{-1}v,l)-x_0(l))\eta\|^2|\xi_V(v)\xi_L(l)|^2\,d\nu(v)d\rho(l)\\
	=&\int_{E_j\times\overline{L_0}}\|(x(v,l)-x_0(l))\eta\|^2|\xi_V(h_jv)\xi_L(l)|^2\,d((h_j^{-1})_*\nu)(v)d\rho(l)\\
	=&\int_{E_j}\left( \int_{\overline{L_0}}\|(x(v,l)-x_0(l))\eta\|^2|\xi_L(l)|^2dl\right)|\xi_V(h_jv)|^2\, d((h_j^{-1})_*\nu)(v)\\
	< &\int_{h_jE_j}\epsilon|\xi_V(v)|^2\,d\nu(v)
	\leq\epsilon\|\xi_V\|^2.
	\end{align*}
	So $1_{h_j E_j}(\sigma_{h_j}(x)-x_0)\to 0$ in SOT, and since $\nu(h_jE_j)\to 1$, also $\sigma_{h_j}(x)-x_0\to 0$ in SOT.
	
	Let $y$ be a WOT cluster point of the set $\{\pi(\gamma_j)x\pi(\gamma_j^{-1})\}$. Then $y\in\mathcal{B}$ because $x\in\mathcal{B}$ and conjugation with $JNJ$ and $J\pi(\Gamma)J$ commutes with conjugation with $\pi(\Gamma)$. Also $yP_1$ is a WOT cluster point of
	\begin{align*}
	\{\pi(\gamma_j)x\pi(\gamma_j^{-1})P_1\}
	=&\{\pi(\gamma_j)(J\pi(\gamma_j)J)(J\pi(\gamma_j^{-1})J)x(J\pi(\gamma_j)J)P_1\}\\
	=&\{\pi(\gamma_j)(J\pi(\gamma_j)J)\sigma_{\gamma_j}(x)P_1\}
	\end{align*}
	Since $\sigma_{\gamma_j}(x)\to x_0$ in SOT, $yP_1$ must then also be a WOT cluster point of \[\{\pi(\gamma_j)(J\pi(\gamma_j)J)x_0P_1\}\subset L^\infty(\overline{L_0})\overline{\otimes} B(L^2M),\]
	so $yP_1\in L^\infty(\overline{L_0})\overline{\otimes} B(L^2M)$.
	$P_1yP_1$ is a WOT cluster point of 
	\begin{align*}
	\{P_1\pi(\gamma_j)x\pi(\gamma_j^{-1})P_1\}
	=\{P_1(J\pi(\gamma_j)J)x(J\pi(\gamma_j^{-1})J)P_1\}
	=\{P_1\sigma_{\gamma}(x)P_1\}.
	\end{align*}
	So again since $\sigma_{\gamma_j}(x)\to x_0$,  $P_1yP_1=P_1x_0P_1$.
\end{proof}
\begin{Prop}\label{M=B}
	If $M$ is not isomorphic to $N\rtimes\Gamma$ with isomorphism extending $\pi$, then $\mathcal{B}=M$, hence $N\subset M$ is coamenable.
\end{Prop}	
\begin{proof}
	Assume that $M$ is not isomorphic to $N\rtimes\Gamma$ with isomorphism extending~$\pi$. Then there is a $\gamma_0\in\Gamma\setminus\{e\}$ and an $n\in N$ such that $c_0:=\tau(\pi(\gamma_0)n)\ne 0$.
	
	Let $x$, $x_0$ and $y$ be as in the above lemma. We want to show that $x$ is a constant function.
	Let $\theta\colon \Gamma\to \Aut(M)$, different as in the proof of Lemma \ref{cond exp on B}, be conjugation by $\pi(\cdot)$ and the induced action $\tilde{\theta}\colon G\curvearrowright L^\infty(G/\Gamma)\overline{\otimes}M$ as in Definition \ref{induced representation}.
	
	$\theta$ is ergodic by assumption, hence $\tilde{\theta}$ is ergodic by Lemma \ref{induced representation ergodic}.
	It is still ergodic when restricted to  $\overline{V_0}\rtimes\mathcal{Z}(R_0)$ because $\overline{V_0}\rtimes\mathcal{Z}(R_0)$ is not compact and hence every $\overline{V_0}\rtimes\mathcal{Z}(R_0)$-invariant vector must also be $G$-invariant by the Howe-Moore property of $G$ \cite[Theorem 5.2]{HoweMoore1979}.
	Now \cite[Lemma 3.2]{PetersonRigidity2014} gives us that for every neighborhood $e\in O\subset G$ and $\Gamma_O=\Gamma\cap O(\overline{V_0}\rtimes\mathcal{Z}(R_0))$ we have
	\begin{align}\label{conv}
	\tau(\pi(\gamma_0)n)=J\tau(\pi(\gamma_0)n)J\in\overline{\conv}^{SOT}\{J\pi(\gamma^{-1})\pi(\gamma_0)n\pi(\gamma)J|\,\gamma\in \Gamma_O\}.
	\end{align}

	We want to show now that $[\sigma_{\gamma_0}^0\otimes c_0, P_1 x_0P_1 ]$ is zero. By Lemma \ref{lem 11.1} $[\sigma_{\gamma_0}^0\otimes c_0, P_1 x_0P_1 ]=[\sigma_{\gamma_0}^0\otimes c_0, P_1 yP_1 ]$.
	The approximation (\ref{conv}) of $\tau(\pi(\gamma_0)n)$ gives for every $O$ an approximation  
	\begin{align*}
	[\sigma_{\gamma_0}^0\otimes c_0, P_1 yP_1 ]
	\overset{SOT}{\sim} \sum_{i=1}^k c_iP_1 \,
	[\sigma_{\gamma_0}^0\otimes J\pi(\gamma_i^{-1}\gamma_0)n\pi( \gamma_i)J, y]\,P_1
	\end{align*}
	with $\gamma_i\in\Gamma_O$ and $\sum_{i=1}^k c_i=1$. 
	
	Now write $\gamma_i=g_ih_i$ where $g_i\in O$ and $h_i\in\overline{V_0}\rtimes\mathcal{Z}(R_0)$.  We have	$\sigma_{h_i}(yP_1 )=yP_1$ and $\sigma_{h_i}(P_1 y)=P_1 y$ since $yP_1, P_1 y \in L^\infty(\overline{L_0})\overline{\otimes} B(L^2M)$ and taking $O$ small enough we get $\sigma_{\gamma_i}(y)P_1 =\sigma_{\gamma_i}(yP_1 )\sim yP_1 $ and $P_1 \sigma_{h_i}(y)=\sigma_{h_i}(P_1y)\sim P_1y$ in SOT.
	Then
	{\allowdisplaybreaks
	\begin{align*}
	&[\sigma_{\gamma_0}^0\otimes c_0, P_1 x_0P_1 ]\\
	\overset{\textup{WOT}}{\sim}& \sum_{i=1}^k c_iP_1 \,
	[\sigma_{\gamma_0}^0\otimes J\pi(\gamma_i^{-1}\gamma_0)n\pi( \gamma_i)J, \sigma_{\gamma_i}(y)]\,P_1\\
	=&\sum_{i=1}^k c_iP_1\,
	[\sigma_{\gamma_0}^0\otimes J\pi(\gamma_i^{-1}\gamma_0\gamma_i)J, \sigma_{\gamma_i}(y)]\,
	(J\alpha_{\gamma_i^{-1}}(n)J)P_1\\
	=& \sum_{i=1}^k c_iP_1 \,
	[\sigma_{\gamma_0}^0\otimes J\pi(\gamma_i^{-1}\gamma_0\gamma_i)J, (\sigma^0_{\gamma_i}\otimes 1)y(\sigma^0_{\gamma_i^{-1}}\otimes 1)]\,
	(J\alpha_{\gamma_i^{-1}}(n)J)P_1\\
	=& \sum_{i=1}^k c_iP_1 (\sigma^0_{\gamma_i}\otimes 1)\,
	[\sigma_{\gamma_i^{-1}\gamma_0\gamma_i}^0\otimes J\pi( \gamma_i^{-1}\gamma_0\gamma_i)J, y](\sigma^0_{\gamma_i^{-1}}\otimes 1)\,
	(J\alpha_{\gamma_i^{-1}}(n)J)P_1\\
	=&0.
	\end{align*}}
	In the second step we used that $y$ and hence also $\sigma_{\gamma_i}(y)$ commutes with $JNJ$. In the last step we used that $[\sigma_{\gamma_i^{-1}\gamma_0\gamma_i}^0\otimes J\pi( \gamma_i^{-1}\gamma_0\gamma_i)J, y]=0$ because $y\in\mathcal{B}$.
	So we found $[\sigma_{\gamma_0}^0\otimes c_0, P_1 x_0P_1 ]=0$ and hence $\sigma_{\gamma_0}(P_1 x_0P_1 )=P_1 x_0P_1$. Since $\tau(\pi(\gamma_0)n)=\tau(\pi(\gamma)\pi(\gamma_0)n\pi(\gamma^{-1}))=\tau(\pi(\gamma\gamma_0\gamma^{-1})\alpha_{\gamma}(n))$, we get $\sigma_{\gamma}(P_1 x_0P_1 )=P_1 x_0P_1 $ for all $\gamma\in \left\langle \left\langle \gamma_0\right\rangle \right\rangle $ in the normal closure of $\gamma_0$. By \cite[Theorem 10.10.]{PetersonRigidity2016} the action of $\left\langle \left\langle \gamma_0\right\rangle \right\rangle$ on $L^\infty(G/P)$ is ergodic, so $P_1 x_0P_1 \in \CC 1\otimes P_{\hat{1}} $.
	Since $x_0$ was arbitrary in the range of $x$, we conclude that $P_1 xP_1 \in L^\infty(\overline{V_0})\otimes P_{\hat{1}} $ and hence $P_1 \mathcal{B}P_1 \subset L^\infty(\overline{V_0})\otimes P_{\hat{1}}$. This means $\mathcal{B}\subset L^\infty(\overline{V_0})\overline{\otimes} B(L^2M)$ because if $x\in\mathcal{B}$ and $a,b\in M$, we have 
	\[\left\langle x\hat{a},\hat{b}\right\rangle 
	=\left\langle (b^*xa)\hat{1},\hat{1}\right\rangle 
	=\left\langle (P_1 b^*xaP_1 )\hat{1},\hat{1}\right\rangle 
	\in L^\infty(\overline{V_0})\] 
	since $b^*xa\in\mathcal{B}$.
	But $\overline{V_0}=G/P_0$ and $G$ is generated by the $P_0$'s \cite[Proposition I.1.2.2]{MargulisDiscrete1991}, so we get
	\[\mathcal{B}= B(L^2M)\cap (J\pi(\Gamma)J\cup JNJ)'=M.\]
	Now $N\subset  M$ is coamenable by Lemma \ref{cond exp on B}.
\end{proof}

\begin{Thm}\label{Thm from Pet machine}
	Let $\Gamma$ be a lattice in a simple real Lie group $G$ which has trivial center and real rank at least $2$. Let $M$ be a finite factor, $N\subset M$ a subfactor and $\pi\colon \Gamma \to \mathcal N_M(N)$ a representation of $\Gamma$ into the normalizer of $N$ such that the action $\Gamma\curvearrowright M$ given by $\alpha_{\gamma}(x)=\pi(\gamma)x\pi(\gamma^{-1})$ is ergodic and $M=(N\cup\pi(\Gamma))''$.
	
	Then $M$ is isomorphic to $N\rtimes \Gamma$ with isomorphism extending $\pi$ or  $[M:N]<\infty$.
\end{Thm}
\begin{proof}
	If $\pi$ does not extend to an isomorphism $M\cong N\rtimes \Gamma$,  the inclusion $N\subset M$ is coamenable by Proposition \ref{M=B}.
   Then by the proof of \cite[Lemma 2.1]{BMO19}	there is a nonzero projection $q\in N'\cap M$ such that $q(N'\cap M)q$ is completely atomic.  Hence the center of $N'\cap M$ is atomic since the conjugation action of $\Gamma$  on it is ergodic and the existence of $q$ implies that it is not diffuse. Since the action is also trace preserving, it is finite. $N'\cap M$ must be of type I since it contains a minimal projection and the action is ergodic, hence it is finite dimensional.
   
	$N\subset M$ is also corigid because $\Gamma$ has property (T) \cite[4.1.7 (ii)]{Popa}.
	So the inclusion is of finite index by Theorem \ref{Popa}.
\end{proof}

In the case where $M=L\Gamma$ and $\pi$ is the left regular representation we have $M\cong N\rtimes \Gamma$ with isomorphism extending $\pi$ if and only if $N=\CC$. Hence we get

\begin{Cor}\label{rigidity left regular}
	Let $\Gamma$ be a lattice in a simple real Lie group $G$ which has trivial center and real rank at least $2$. Let $N\subset L\Gamma$ a subfactor which is normalized by the natural copy of $\Gamma$ in $L\Gamma$. Then $N=\CC$ or $[L\Gamma:N]<\infty$.
\end{Cor}
Corollary \ref{rigidity left regular} can also be obtained without Theorem \ref{Thm from Pet machine} from \cite[Theorem 3.15 1)]{CD19} and \cite[Section 4.3]{Doktorarbeit}.
In \cite{CD19} similar results are proven for groups with positive $\ell^2$-Betti numbers and acylindrically hyperbolic groups.
S.\@ Das brought to our attention that S.\@ Popa observed that in the above situation the index is an integer as is always true for regular subfactors with finite index, see also \cite[Theorem 4.5]{C-S06}.

\section{On deducing the theorem from character rigidity}
It seems natural to try to deduce Theorem \ref{Thm from Pet machine} from character rigidity \cite[Theorem C]{PetersonRigidity2016} by applying it to the character $\gamma\mapsto \|E_N(\pi(\gamma))\|_2^2$. The following connection is easy to deduce:

\begin{Lemma}\label{extremal phi}
	Let $N\subset M$ and $\pi$ be as above and let
$\varphi(\gamma)= \|E_N(\pi(\gamma))\|_2^2$.
If the GNS construction of  this character generates a finite dimensional von Neumann algebra, then $[M:N]<\infty$. If it is the regular character, then $M=N\rtimes\Gamma$.  
\end{Lemma}
\begin{proof}
If $\varphi$ is the regular character, \[
E_N(\pi(\gamma))=\begin{cases}
0&\gamma\ne e\\
1&\gamma=e
\end{cases},\]
hence $M=N\rtimes\Gamma$. 

We now describe the GNS construction of $\varphi$.
Let $e_N\in B(L^2M)$ be the orthogonal projection onto $N$  and $\langle M,N\rangle:=\{\sum_{k=1}^n x_ke_Ny_k|\,n\in \NN, x_k,y_k\in M\}''\subset B(L^2M)$ the basic construction with semifinite trace given by $Tr(xe_Ny)=\tau(xy)$. From this we get a Hilbert space $\mathcal{H}:=L^2(\langle M,N\rangle, Tr)$ as the  completion of the finite elements of $\langle M,N\rangle$ with the norm $\|x\|_2=Tr(x^*x)^{\frac{1}{2}}$. Define a unitary representation $\theta\colon \Gamma\to U(\mathcal{H})$ by
\[\theta_\gamma(xe_Ny)=\pi(\gamma)xe_Ny\pi(\gamma)^*.\]
Then, using $e_Nxe_N=E_N(x)e_N$ for all $x\in M$, we get
\begin{align*}
\left<\theta_\gamma(e_N),e_N\right>&=\left<\pi(\gamma)e_N\pi(\gamma)^*,e_N\right>=Tr(e_N\pi(\gamma)e_N\pi(\gamma)^*e_N)=Tr(E_N(\pi(\gamma))e_N\pi(\gamma)^*e_N)\\
&=Tr(E_N(\pi(\gamma))e_NE_N(\pi(\gamma)^*))=\tau_M(E_N(\pi(\gamma))E_N(\pi(\gamma))^*)=\tau(\gamma),
\end{align*}
hence $\mathcal{H}_\varphi:= \overline{ \spa} \{\sigma_\gamma(e_N)|\gamma\in\Gamma\}\subset \mathcal{H}, e_N$, and $\sigma$ form a GNS triple for $\varphi$.

Assume $\sigma(\Gamma)''\subset B(\mathcal{H}_\varphi)$ is finite dimensional. Say it is generated as a vector space by $\sigma_{\gamma_1},\dots,\sigma_{\gamma_n}$. Then for all $x,y\in N$
\[\pi(\gamma)xe_Ny\pi(\gamma')=\pi(\gamma)e_Nxy\pi(\gamma)^*\pi(\gamma\gamma')
=\theta_\gamma(e_Nxy)\pi(\gamma\gamma')
=\sum_{i=1}^n c_i\pi(\gamma_i)e_Nxy\pi(\gamma_i)^*\pi(\gamma\gamma')\]
for some $c_i\in\CC$.
So since $M$ is generated by $N$ and $\pi(\Gamma)$, $\langle M,N\rangle$ is generated over $M$ by $\pi(\gamma_1)e_N, \dots,\pi(\gamma_n)e_N$. Hence $[M:N]=[\langle M,N\rangle: M]<\infty$. 
\end{proof}

In particular, if we knew that $\varphi$ was extremal or could reduce the situation to the extremal case, the theorem would follow. However, things are more complicated.
\begin{Def}
	A trace-preserving action on a finite von Neumann algebra $\sigma\colon\Gamma\curvearrowright M$  that leaves a von Neumann subalgebra $N\subset M$ invariant is  called \textit{weakly mixing relative to $N$} if for any finite set $F\subset M$ with $E_N(x)=0$ for all $x\in F$   there exist $\gamma_n\in\Gamma$ such that for all $\eta,\eta'\in F$, $\|E_N(\eta^*\sigma_{\gamma_n}(\eta'))\|_2\to 0$ for $n\to\infty$.
\end{Def}
\begin{Lemma}
	The following are equivalent.
	\begin{enumerate}
		\item $\varphi$ is extremal.
		\item $\sigma_\gamma(x)=\pi(\gamma)x\pi(\gamma)^*$ is weakly mixing relative to $N$.
		\item $M=N\rtimes\Gamma$.
	\end{enumerate}
\end{Lemma}
\begin{proof}
$\varphi$ is extremal if and only if $\theta(\Gamma)''\subset B(L^2(\left<M,N\right>,Tr))$ is a factor. We have a dense inclusion $\theta(\Gamma)''\subset L^2(\theta(\Gamma)'')\cong\mathcal{H}_\varphi$ sending $\theta(\gamma)\in\theta(\Gamma)''$ to $\pi(\gamma)e_N\pi(\gamma)^*\in\mathcal{H}_\varphi$ such that the conjugation action on $\theta(\Gamma)''$ translates to the action $\sigma^N$ on $\mathcal{H}_\varphi$ given by
\[\sigma^N_\gamma(xe_Ny)= \pi(\gamma)x\pi(\gamma)^*e_N\pi(\gamma)y\pi(\gamma)^*.\]
So $\varphi$ is extremal iff every $\sigma^N$-invariant vector in $\mathcal{H}_\varphi$ is trivial. Since $N^\sigma=\CC$ this is equivalent to 
 \begin{align}\label{iii'}
 (\mathcal{H}_\varphi)^{\sigma^N}\subset L^2(Ne_N). 
  \end{align}
 Condition (\ref{iii'}) is equivalent to $\sigma$ being weakly mixing because it is implied by $(iii)$ in \cite[Lemma 2.10]{Popa2007} and implies $(i)$ with the same proof as for $(iii)\Rightarrow (i)$. Hence i) and ii) are equivalent.

If $\varphi$ is extremal, then $M=N\rtimes \Gamma$ or $[M:N]$ by character rigidity and Lemma \ref{extremal phi}. In the second case,	since $N$ is a factor, $L^2(M)$ has a finite orthogonal basis $\eta_1,\dots,\eta_k$ over $N$ with 
	$\|\eta\|_2=\sum_{i=1}^k\|E_N(\eta_i^*\eta)\|_2$ for  and all $\eta\in L^2(M)$ (see \cite[8.4-8.6]{ADPopa}).
	Then for all $\gamma\in\Gamma$, 
	\[\sum_{i=1}^k\|E_N(\eta_i^*\sigma_\gamma(\eta)\|_2
	=\|\sigma_\gamma(\eta)\|_2=\|\eta\|_2,	\]
	hence $\|E_N(\eta_i^*\sigma_{\gamma_n}(\eta))\|_2$ cannot go to zero for some $\gamma_n$ and all $i=1,\dots,k$.
\end{proof}
\begin{Ex}
	If $\Lambda\triangleleft\Gamma$ is a normal subgroup, $\pi\colon\Gamma\to L\Gamma=M$ the left regular representation and $N=L\Lambda$, then $\varphi$ is the regular character on $\Gamma/\Lambda$ which is not extremal if the index is finite.
	The decomposition into extremal characters corresponds to the decomposition of the left regular representation of $\Gamma/\Lambda$ into irreducible representations, which doesn't seem to be nicely reflected in the situation of $\Gamma\curvearrowright N\subset M$.
\end{Ex}
\phantomsection
\addcontentsline{toc}{chapter}{Bibliography}
\bibliographystyle{alpha}

\end{document}